\documentclass[11pt,oneside]{amsart}
\usepackage{geometry} 
\usepackage{amsfonts}
\usepackage{epsfig}
\usepackage{url}
\usepackage{amssymb,latexsym,amsmath,amsthm,verbatim}
\usepackage{graphicx,epsfig,epstopdf,amssymb,color}
\usepackage{bm}
\newcommand{\veps}{\varepsilon}
\newcommand{\R}{\mathbb{R}} 
\newcommand{\C}{\mathbb{C}}
\newcommand{\K}{\mathbb{K}}
\newcommand{\OO}{\mathcal{O}}

\newcommand{\ba}{\begin{array}}
\newcommand{\ea}{\end{array}}
\newcommand{\tPsi}{\widetilde{\Psi}}
\newcommand{\tA}{\widetilde{A}}
\newcommand{\EE}{{\bf E}}
\newcommand{\eps}{\varepsilon}

\newtheorem{theorem}{Theorem}[section]
\newtheorem{lemma}[theorem]{Lemma}

\newtheorem{corollary}[theorem]{Corollary}
\newtheorem{remark}[theorem]{Remark}

\begin{document}
\title{Justification of the Nonlinear Schr\"odinger approximation  
for a quasilinear Klein-Gordon equation}
\author{Wolf-Patrick D\"ull$^1$}\thanks{$^1$IADM, 
Universit\"at Stuttgart, Pfaffenwaldring 57, 70569 Stuttgart, Germany\\
(duell@mathematik.uni-stuttgart.de)}
\date{\today}
\maketitle
\begin{abstract}
We consider a nonlinear Klein-Gordon equation with a quasilinear quadratic term.
The Nonlinear Schr\"odinger (NLS) equation can be derived as a formal approximation equation describing the evolution of the envelopes of slowly modulated  spatially and temporarily oscillating wave packet-like solutions to the quasilinear Klein-Gordon equation. 
It is the purpose of this paper to present a method which allows one to prove error estimates in Sobolev norms between exact solutions of the quasilinear Klein-Gordon equation and the formal approximation obtained via the NLS equation.
The paper contains the first validity proof of the NLS approximation 
of a nonlinear hyperbolic equation with a quasilinear quadratic term by error estimates in Sobolev spaces.
We expect that the method developed 
in the present paper will allow an answer to the relevant question of the validity of the NLS approximation for other quasilinear hyperbolic systems.
\end{abstract}


\section{Introduction and Result}

\begin{figure}[htbp]
\epsfig{file=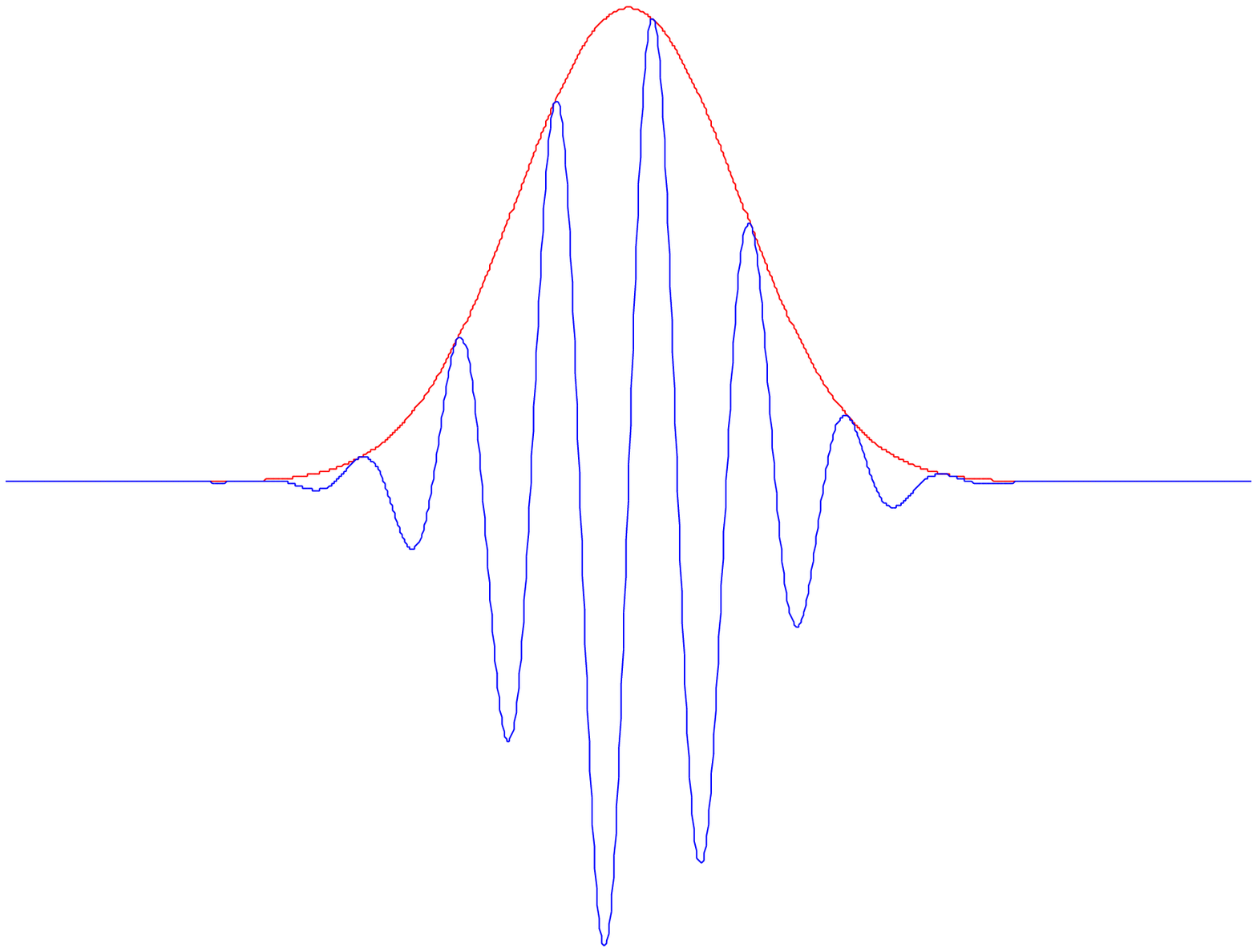,width=15cm,height=6.2cm,angle=0}
\vspace*{-6.3cm}

\hspace{1.9cm}
\epsfig{file=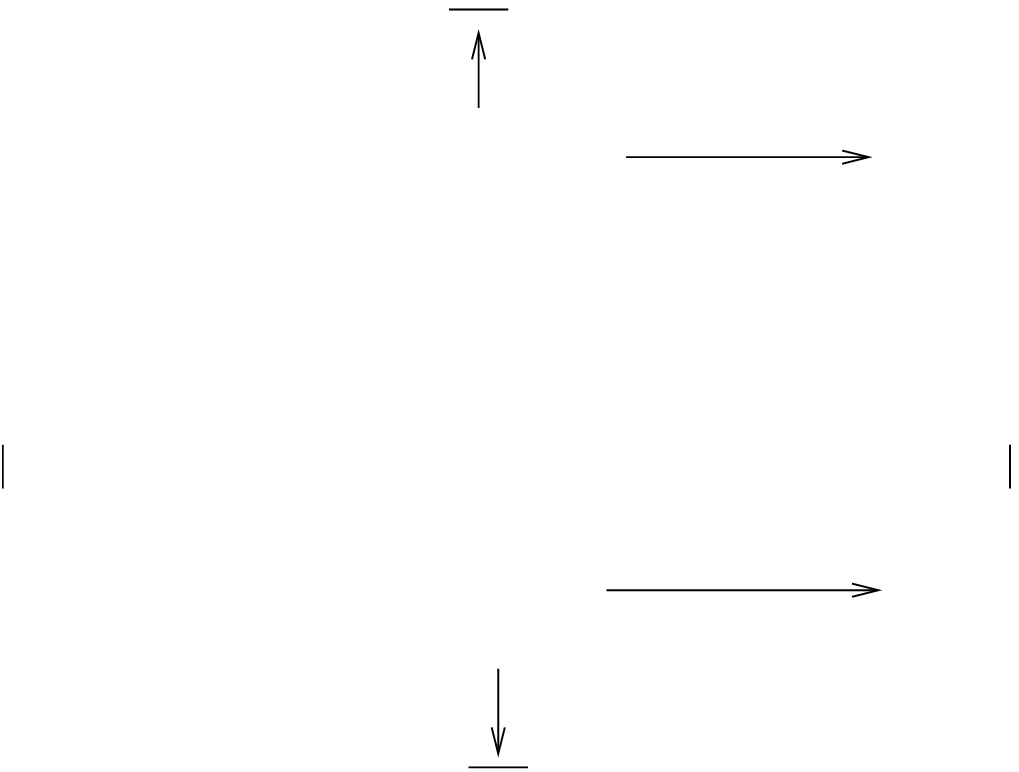,width=10.7cm,height=6.5cm,angle=0}

\vspace*{-5.34cm}
\hspace*{9.7cm}$c_\mathrm{g}$  \vspace*{3.2cm}

\hspace*{9.7cm} $c_\mathrm{p}$

\vspace*{-2.5cm}
\hspace{2.2cm}\hspace{5.0cm}$\varepsilon$
\vspace*{0.3cm}

\hspace{2.65cm}$1/\varepsilon$
\hspace{2.7cm}
\vspace*{2.6cm}

\caption{{\small The envelope (advancing with the group velocity 
$ c_g $) of the oscillating wave packet 
(advancing with the phase velocity 
$ c_p= \omega_0/k_0  $) is described by the 
amplitude $ A $ which solves  the NLS equation \eqref{NLS}.} \label{fig1} }
\end{figure}

The Nonlinear Schr\"odinger (NLS) equation plays an important role in describing  approximately  
slow modulations in time and space of an underlying spatially and temporarily oscillating wave packet in a more complicated hyperbolic system, such as Maxwell's equations for modeling nonlinear optics or the equations describing surface water waves, see, for example, \cite{AS81}.
In this paper, we study the NLS approximation of the quasilinear Klein-Gordon equation 
\begin{equation} \label{qweq}
\partial_t^2 u = \partial_x^2 u -u + \partial_x^2(u^2 )\,, 
\end{equation}
with $ x,t \in \R $, and $ u(x,t) \in \R $.  We make the ansatz $u=\veps \Psi_{NLS} + \OO(\veps^2)$, with
\begin{equation} \label{ansatz}
\veps \Psi_{NLS}(x,t) 
= \veps A(\eps
(x-c_g t),\eps^2t) e^{i( k_0 x - \omega_0 t)} + \mathrm{c.c.} \,.
\end{equation}
Here $0 < \eps \ll 1$ is a small perturbation parameter,
$ \omega_0 > 0$ the basic temporal 
wave number associated to the basic spatial wave number $ k_0 > 0$ of the underlying carrier wave $ e^{i(k_0 x - \omega_0 t)}$,
$c_g$ the group velocity, $A$ the complex-valued amplitude, and c.c. the complex conjugate. 
With the help of (\ref{ansatz}) we describe slow spatial and temporal modulations of the underlying carrier wave.
Inserting the above ansatz into \eqref{qweq} we find that $A$ satisfies at leading order in $\eps$ the NLS equation
\begin{equation} \label{NLS}
\partial_T A  = i \nu_1 \partial_{X}^2 A + i  \nu_2 A|A|^2\,, 
\end{equation}
where $X = \eps(x-c_g t)$, $ T = \eps^2 t $, and    $\nu_j= \nu_j(k_0) \in \R$. $ T $  is the slow time scale and $ X$
is the slow spatial scale, that means,
the time scale of the modulations is  
$\OO({1/\eps^2})$ and the spatial scale of the modulations
is $\OO({1/\eps})$. See Figure \ref{fig1}.
The basic spatial wave number $k= k_0$ and the basic temporal
wave number $\omega = \omega_0$ 
are related via  the linear dispersion relation
of the quasilinear Klein-Gordon equation (\ref{qweq}), namely 
\begin{equation}
\label{lindis}
\omega^2(k) - (1+ k^2) =0\,,
\end{equation}
where we choose the branch of solutions
\begin{equation}
\label{branchomega}
 \omega(k) :=   \ {\rm
 sign}(k)\sqrt{1+ k^2}\,.
\end{equation}
Then the group velocity $ c_g $ of the wave packet
is given by $ c_g = \partial_k  \omega|_{k=k_0}$.
Our ansatz leads to waves moving to the right. To obtain waves moving
to the left, $-\omega_0$ and $c_g$ have to be replaced by  $\omega_0$ and $-c_g$.

\vspace*{0.35cm}

It is the goal of the present paper to prove 
\begin{theorem}
\label{mainresult} 
Fix $s_A \geq 6$.  Then for all $k_0 > 0$ 
and for all $C_1,T_0 > 0$ there exist $C_2 > 0$, $\eps_0 > 0$ 
such that for all solutions $A \in
C([0,T_0],H^{s_A}(\R,\C))$ of the NLS equation (\ref{NLS})
with 
$$ 
\sup_{T \in [0,T_0]} \| A(\cdot,T) \|_{H^{s_A}(\R,\C)} \leq C_1
$$  
the following holds.
For all $\eps \in (0,\eps_0)$
there are solutions 
$$
u \in
C([0,T_0/\eps^2], H^{s_A}(\R,\R))
$$
of the quasilinear Klein-Gordon equation \eqref{qweq} which satisfy
  $$\sup_{t\in[0,T_0/\eps^2]} \| u(\cdot,t) -
    \eps \Psi_{NLS}(\cdot,t)\|_{H^{s_A}(\R,\R)} 
\le C_2 \eps^{3/2}.$$
\end{theorem}

The error of order $\OO(\eps^{3/2})$ is small compared with the solution 
$u$ and the approximation $\eps \Psi_{NLS}$, which are both of order
$\OO(\eps)$ in $ L^{\infty} $ such that the dynamics 
of the NLS equation can be found in the quasilinear Klein-Gordon
equation \eqref{qweq}, too. The NLS equation is a completely integrable
Hamiltonian system, which can be solved 
explicitly with the help of some inverse scattering scheme, see, for example, \cite{AS81}.

It should be noted that the smoothness in our error bound is equal to the assumed smoothness of the amplitude. This can be achieved by using a modified approximation which has compact support in Fourier space but differs only slightly from $\eps \Psi_{NLS}$. Such an approximation can be constructed because the Fourier transform of $\eps \Psi_{NLS}$ is sufficiently strongly concentrated around the wave numbers $\pm k_0$.

Like in many other proofs of error estimates in the literature we will assume in our proof of Theorem \ref{mainresult} that $s_A$ is an integer in order to simplify the analysis by using Leibniz's rule,  but our proof can be generalized to be valid for all $s_A \geq 6$.

We remark that such an approximation theorem like Theorem \ref{mainresult} should not be taken for granted. There are various counterexamples,
where approximation equations derived by reasonable formal arguments
make wrong predictions about the dynamics of the 
original systems, see, for example, \cite{Schn05,SSZ15}. For an introduction into theory and applications of the NLS approximation we refer to \cite{Schn11OWbuch}.

\medskip

In order to explain our method to prove Theorem \ref{mainresult} and the relevance of this method we
consider the more general abstract evolutionary
problem
\begin{equation} \label{abstrsyst}
  \partial_t V = \Lambda V + B(V,V)\,,
\end{equation}
with $V=V(x,t) \in \R^2$, $x,t \in \R$,
$\Lambda$ a linear operator whose symbol is a diagonal matrix of the form 
\begin{equation} \label{linsymb}
\widehat{\Lambda}(k) = \mathrm{diag}\,(-i\lambda(k),\, i\lambda(k))\,,
\end{equation}
where $k\in \R$ and $\lambda$ is a piecewise smooth real-valued odd function,
and $B$ a symmetric bilinear operator. 
Writing the quasilinear Klein-Gordon equation \eqref{qweq} as a first-order system and diagonalizing the linear part of the resulting system we obtain a special case of \eqref{abstrsyst}, 
where $\lambda =\omega$, with $\omega$ defined by \eqref{branchomega}.

The NLS equation \eqref{NLS} can be derived as a formal approximation equation with the help of the ansatz $V=\eps \tilde{\Psi}$, where
\begin{align} \label{approx1}
\eps \tilde{\Psi}(x,t) =\;& \eps \left(\begin{array}{c}
\Psi_{NLS}(x,t)\\
0
\\ \end{array}\right)
+ \eps^2 \left(\begin{array}{c}
 {\tA}_{01} (\eps (x-c_gt),\eps^2t)\\
{\tA}_{02} (\eps (x-c_gt),\eps^2t)
\\ \end{array}\right)\\[2mm] & 
+ \eps^2 \left( \left(\begin{array}{c}
 {\tA}_{21} (\eps (x-c_gt),\eps^2t)\\
{\tA}_{22} (\eps (x-c_gt),\eps^2t)
\\ \end{array}\right) e^{2i( k_0 x - \omega_0 t)} + \mathrm{c.c.} \right), \nonumber
\end{align}
with $\Psi_{NLS}$ as in \eqref{ansatz}, where $\omega_0 = \lambda(k_0)$ and $c_g=\partial_k \lambda(k_0)$, real-valued functions ${\tA}_{01}, {\tA}_{02}$, and complex-valued functions ${\tA}_{21}, {\tA}_{22}$. Inserting this ansatz into \eqref{abstrsyst} and equating the coefficients in front of the $\veps^m  e^{ji( k_0 x - \omega_0 t)}$ for $m \in \{1,2,3\}$ and  $j \in \{0,1,2\}$ to zero yields the NLS equation \eqref{NLS} if $\lambda$ satisfies $\lim\limits_{k \to 0^{\pm}} \lambda(k) \neq 0$ or $\lim\limits_{k \to 0^{\pm}} \partial_k \lambda(k) \neq c_g$ as well as $\pm \lambda(2k_0)\neq 2 \lambda(k_0)$, which is true for $\lambda=\omega$.

It is possible to modify  $\eps \tilde{\Psi}$ to make it an even more accurate approximation. Indeed, for all $\gamma > 0$ there 
exists a function $\Psi$ such that $\Psi-\tilde{\Psi} \to  0$ for $\veps \to 0$ and
\begin{equation} \label{modi}
\mbox{Res}(\veps \Psi) := - \partial_t (\veps \Psi) + \Lambda (\veps \Psi) 
+ B(\veps \Psi, \veps \Psi) = \OO(\eps^\gamma)\,. 
\end{equation}

In order to prove Theorem \ref{mainresult} we have to estimate the error
 \begin{equation}
\eps^\beta R := V - \eps {\Psi}
\end{equation}
for all $t \in [0,T_0/\eps^2]$ to be of order $\OO(\eps^\beta)$ for a $ \beta > 1 $, 
that means, we have to prove that
$R$ is of order $\OO(1)$ for all $t \in [0,T_0/\eps^2]$.
The error $R$ satisfies the equation
\begin{equation}
  \partial_tR = \Lambda R + 2\eps^\alpha B({\Psi},R) + \eps^\beta B(R,R) +
    \eps^{-\beta} \mbox{Res}(\eps {\Psi})\,,
\end{equation}
with $\alpha =1$.

Since our linear operator $\Lambda$ generates a uniformly bounded strongly continuous semigroup, we would be done 
if we had a) $\alpha \geq 2$, b) $\beta > 2$, and c) $
\eps^{-\beta}{\rm 
Res}(\eps{\Psi}) = \OO(\eps^2)$.
The result then would follow by a rescaling of time, $T = \eps^2t$, and an application
of Gronwall's inequality, see, for example, \cite{KSM92}.  However, we have $\alpha = 1$.
We can still make  $\gamma$ in (\ref{modi})  arbitrary large by constructing our approximate solution $\veps \Psi$ as described below,  and in particular, 
strictly bigger than
$4$.  Consequently, we can choose $\beta > 2$ and so the points b) and c) are satisfied easily.
Hence, the difficulty is to control the quadratic term $2\eps B({\Psi},R)$.

Semilinear quadratic terms can be eliminated with the help of a so-called normal-form transform 
\begin{equation} \label{inft} 
\tilde{R} := R + \eps
N({\Psi},R)\,,
\end{equation}
 where $N$ is an appropriately constructed bilinear mapping, if the so-called non-resonance condition 
\begin{equation} \label{inonres}
\inf_{k\in \R,\atop  j_1, j_2 \in \{\pm1\}}\, | {-j_1\lambda(k)-\lambda(k_0)+j_2\lambda(k-k_0)}| \geq  C > 0
\end{equation}
is satisfied. This normal-form transform is invertible and the new error function $\tilde{R}$ satisfies an evolution equation of the form 
\begin{equation} \label{newerr}
\partial_t \tilde{R} = \Lambda \tilde{R} + \eps^2 g({\Psi},\tilde{R}) +
    \eps^{-\beta} \mbox{Res}(\eps {\Psi})\,,
\end{equation}
where $g({\Psi},\tilde{R})$ is a semilinear term of order $\OO(1)$. Therefore, to this equation, Gronwall's inequality can be applied to bound $\tilde{R}$ and hence $R$ for all $t \in [0,T_0/\eps^2]$.

The strategy of using normal-form transforms to eliminate semilinear quadratic terms in hyperbolic systems was introduced in \cite{Sh85}. In the context
of justifying NLS approximations, it was first applied in \cite{Kal88} and was further developed in \cite{Schn98Nodea, Schn98JDE}, where also first attempts were made to weaken the non-resonance condition \eqref{inonres}. Proving the validity of 2D NLS approximations with the help of a normal-form transform was addressed in \cite{DHSZ15}.

However, in the case of equation \eqref{qweq} the quadratic term $2\eps B({\Psi},R)$ is quasilinear and loses one derivative, that means, $ R \mapsto 2\eps B({\Psi},R)$ maps $H^{m+1}(\R,\C)$ into $H^{m}(\R,\C)$ or $C^{n+1}(\R,\C)$ into $C^{n}(\R,\C)$. Since the non-resonance condition \eqref{inonres} is satisfied for $\lambda=\omega$, it is still possible to construct a normal-form transform of the form \eqref{inft} to eliminate this term. But the normal-form transform also loses one derivative. At least, the term $2\eps B({\Psi},R)$ being generated by equation \eqref{qweq} has a structure that nevertheless allows one to invert the normal-form transform. Hence, the new error function $\tilde{R}$ still satisfies an evolution equation of the form \eqref{newerr}, but
the loss of one derivative caused by the normal-form transform implies that the term $g({\Psi},\tilde{R})$ loses even two derivatives.

Consequently, the variation of constants formula and Gronwall's inequality cannot be applied to bound $\tilde{R}$ - in this situation not because of missing powers of $\veps$ but due to regularity problems, since the semigroup generated by $\Lambda$ is not smoothing. For the same reasons it does not work either to bound $\tilde{R}$ by deriving energy inequalities for $\tilde{R}$ and applying Gronwall's inequality to them.

Hence, the validity of the NLS approximation for hyperbolic
systems with quasilinear quadratic terms is a highly nontrivial problem, which
has been remained unsolved in general for more than four decades.
Only for some examples of quasilinear systems the NLS approximation 
has been justified so far. All these examples avoid the major difficulty captured by equation \eqref{qweq}. 

The first and very general NLS approximation theorem for quasilinear dispersive
wave systems was proven in \cite{Kal88}. However, the occurrence of quasilinear quadratic 
terms was excluded 
explicitly.
Another example are dispersive wave systems where the right-hand sides lose only half a derivative. The 2D water wave problem without surface tension and finite depth in Lagrangian coordinates falls into this class.
 In this case
the elimination of the quadratic terms is possible with the help of normal-form transforms. The right-hand sides of 
the transformed systems then lose one derivative and can be handled with the help of  the Cauchy-Kowalevskaya theorem \cite{SW10,DSW12}.
Furthermore, the NLS approximation was justified for the 2D and 3D water wave problem without surface tension and infinite depth \cite{TW11,T14} by finding a different 
transform adapted to the special structure of that problem. 
Similarly, for 
the quasilinear Korteweg-de Vries
equation the result can be obtained by simply applying a
Miura transform \cite{Schn11}.
In \cite{CDS15}, the NLS approximation of time oscillatory long waves for equations with quasilinear quadratic terms was proven for analytic data without using a normal-form transform.
Moreover, another approach to address the problem of the validity of the NLS approximation can be found in \cite{MN13}.
Finally, some numerical evidence that the NLS approximation is also valid for quasilinear equations like equation \eqref{qweq} was given in \cite{CS11}. 

The present paper contains the first validity proof of the NLS approximation 
of a nonlinear hyperbolic equation with a quasilinear quadratic term in Sobolev spaces.
We expect that the method developed 
here will allow an answer to the relevant question of the validity of the NLS approximation for many other quasilinear hyperbolic systems.
Our method of proof is as follows.

Instead of performing the normal-form transform \eqref{inft} we only use the term $\eps N(\Psi,R)$ to define the energy 
\begin{equation} \label{ien}
\mathcal{E}_{s} = \sum_{\ell=0}^s  \Big( \frac{1}{2}   \|\partial_x^{\ell} R\|_{L^2}^2  +  \veps \int_{\R}\partial_x^{\ell} R \bullet \partial_x^{\ell} N(\Psi,R)\,dx \Big)\,,  
\end{equation}
where $s= s_A \geq 6$. Since $\|\tilde{R}\|_{H^s}^2$ differs from  $\mathcal{E}_{s}$ only by terms of order $\mathcal{O}(\veps^2)$, the evolution equations of $\mathcal{E}_{s}$ and $\|\tilde{R}\|_{H^s}^2$ share the property that their right-hand sides are of order $\mathcal{O}(\veps^2)$.

This strategy was already used in \cite{D12} as an ingredient to simplify the proof of error estimates 
compared with the alternative proofs in \cite{SW00, SW02}. To overcome regularity problems caused by quasilinear quadratic terms, this strategy was first used in \cite{HITW13} and was further developed in \cite{HIT14, IT14} to apply it to the water wave problem with infinite depth. In these three papers, structural properties of the Hilbert transform help to construct and estimate the energy.  

In the case of equation \eqref{qweq} we can show that $N(\Psi,R)$ can be split into a term of the form $\mathrm{diag}\, (f_1(\Psi),\, f_2(\Psi))\, \partial_x R$ and terms which do not lose regularity such that partial integration yields the
equivalence of $\sqrt{\mathcal{E}_{s}}$ and $\|R\|_{H^s}$ for sufficiently small $\veps$. Consequently, the right-hand side of the evolution equation of $\mathcal{E}_{s}$ can be written as a sum of integral terms containing at most one factor $\partial_x^{s+1} R$ and not two. 

Therefore, the structure of $\Lambda$ and the properties of $\omega$ allow us to construct a modified energy 
\begin{equation} \label{eneq0}
\tilde{\mathcal{E}_{s}} = \mathcal{E}_{s} +\veps^2 h\,,
\end{equation}
where $h=\mathcal{O}(\|R\|_{H^s}^2)$ as long as $\|R\|_{H^s} = \OO(1)$ and 
 $\veps^2 \partial_t h$ eliminates all the 
integral terms on the right-hand side of the evolution equation of $\mathcal{E}_{s}$ with a factor $\partial_x^{s+1} R$.
Consequently, we obtain
\begin{equation}
\label{eneq1}
\partial_t \tilde{\mathcal{E}_s} \lesssim \veps^{2} (\tilde{\mathcal{E}_s}+1)
\end{equation}
as long as $\|R\|_{H^s} = \OO(1)$ such that Gronwall's inequality yields the $\mathcal{O}(1)$-boundedness of $\tilde{\mathcal{E}_s}$ and hence of $R$ for all $t\in[0,T_0/\veps^{2}]$. 

We finish the discussion of our method of proof with a comment on the chosen regularity of the quadratic nonlinearity. The term $\partial_x^2(u^2 )$ is the  
least regular quadratic term which can be added to the linear part of the right-hand side of equation \eqref{qweq} so that the method developed in the present paper allows us to construct an energy $\tilde{\mathcal{E}_s}$  satisfying \eqref{eneq0}-\eqref{eneq1} as long as $\|R\|_{H^s} = \OO(1)$. It is a question of ongoing research if this method can be generalized to be applicable to less regular nonlinearities.

\medskip

The plan of the paper is as follows. In Section \ref{sec:appr} we rewrite the quasilinear Klein-Gordon equation \eqref{qweq} to make it a special case of the evolutionary system \eqref{abstrsyst} and derive the NLS approximation. In Section \ref{sec:err}
we present the error equations, construct our energy and perform the error estimates to prove Theorem \ref{mainresult}.

In a forthcoming paper, we intend to combine the methods of the present paper with the methods from \cite{D12, DS06} to generalize the approximation theorem for the NLS approximation from \cite{DS06} to the case of the water wave problem with surface tension.  

\medskip
{\bf Notation}. 
We denote the 
Fourier transform of a function $u \in L^2(\R,\K)$, with $\K=\R$ or $\K=\C$ by
$$\hat{u}(k) = \frac{1}{2\pi} \int_{\R} u(x) e^{-ikx} dx. $$
Let $H^{s}(\R,\K)$ be
the space of functions mapping from $\R$ into $\K$
for which
the norm
$$ \| u \|_{H^{s}(\R,\K)} = \left(\int_{\R} |\hat{u}(k)|^2 (1+|k|^2)^{s} 
dk \right)^{1/2} $$ 
is finite. We also write $L^2$ and $H^{s}$ instead of $L^2(\R,\R)$ and $H^{s}(\R,\R)$.
Moreover, we use the space
$ L^p(m)(\R,\K) $ defined by $ u \in L^p(m)(\R,\K) \Leftrightarrow u \sigma^m \in L^p(\R,\K)$, where
$ \sigma(x) = (1+x^2)^{1/2}$.

Furthermore, we write $A
\lesssim B$, if $A \leq C B$ for a constant $C>0$, and $A = \mathcal{O}(B)$, if 
$|A| \lesssim B$. 

\section{The Derivation of the NLS Approximation}
\label{sec:appr}
We rewrite the quasilinear Klein-Gordon equation
\begin{equation*} 
\partial_t^2 u = \partial_x^2 u -u + \partial_x^2(u^2 ) 
\end{equation*}
as a first-order system
\begin{align}
\partial_t u & =   \sqrt{1-\partial_{x}^2}\, Hv\,,\\
\partial_t v & =  \sqrt{1-\partial_{x}^2}\, Hu - \partial_{x}^2 (1-\partial_{x}^2)^{-1/2} H(u^2)\,,
\end{align}
where $H$ denotes the Hilbert transform. In Fourier space we have
\begin{align}
\partial_t \widehat{u}(k,t) & =    -i \omega(k) \widehat{v}(k,t)\,,\\
\partial_t \widehat{v}(k,t) & =   -i \omega(k) \widehat{u}(k,t) - i \rho(k) (\widehat{u} \ast \widehat{u})(k,t) \,,
\end{align}
where
\begin{equation} \label{omegasym}
\omega(k) = \mathrm{sign}(k) \sqrt{1+k^2}
\end{equation}
and
\begin{equation} \label{rhosym}
\rho(k) = \mathrm{sign}(k) \frac{k^2}{\sqrt{1+k^2}}\,.  
\end{equation}
We diagonalize this system by 
\begin{equation}  
 \left(
\begin{array}{c} \widehat{u} \\ \widehat{v}
\end{array}
\right) = 
\left(
\begin{array}{cc}
  1 &  1   \\
 1 & -1    
\end{array}
\right)
\left(
\begin{array}{c} \widehat{u}_{-1} \\ \widehat{u}_{1} 
\end{array}
\right).
 \end{equation}
Then we obtain
\begin{equation}  
 \left(
\begin{array}{c} \widehat{u}_{-1} \\ \widehat{u}_{1}
\end{array}
\right) = 
\frac{1}{2}\left(
\begin{array}{cc}
  1 &  1   \\
 1 & -1    
\end{array}
\right)
\left(
\begin{array}{c} \widehat{u} \\ \widehat{v}
\end{array}
\right)
 \end{equation}
and
\begin{align} \label{diag1}
\partial_t \widehat{u}_{-1} (k,t) & =    -i \omega(k) \widehat{u}_{-1}(k,t)
 - \frac{1}{2} i \rho(k) (\widehat{u}_{-1}+\widehat{u}_1)^{\ast 2}(k,t)\,,\\ \label{diag2}
\partial_t \widehat{u}_1(k,t) & =   i \omega(k) \widehat{u}_1(k,t) + \frac{1}{2} i \rho(k) (\widehat{u}_{-1} + \widehat{u}_1)^{\ast 2}(k,t) \,.
\end{align}
\\
In order to derive the NLS equation as an approximation equation for system (\ref{diag1})-(\ref{diag2}) we make the ansatz
\begin{equation} \label{ans1c}
\left( \ba{c}
{u}_{-1}\\ {u}_1 \ea \right) = \veps \tPsi =
\eps {\tPsi}_{1} + \eps {\tPsi}_{-1} +  \eps^2 {\tPsi}_0 + \eps^2{\tPsi}_2 
+ \eps^2{\tPsi}_{-2}\,,
\end{equation}
with
\begin{align*}
\eps {\tPsi}_{\pm 1} 
& =  \eps {\tA}_{\pm 1} (\eps
(x-c_g t),\eps^2 t)
\,\EE^{\pm 1} \!\left( \ba{c} 
1 \\ 0  \ea \right),\\[2mm]
\eps^2 \tPsi_0 
& =  \left(\begin{array}{c}
\eps^2 {\tA}_{01} (\eps (x-c_gt),\eps^2t)\\
\eps^2 {\tA}_{02} (\eps (x-c_gt),\eps^2t)
\\ \end{array}\right),
\end{align*}
\begin{align*}
\eps^2 {\tPsi}_{\pm2} 
& =  \left(\begin{array}{c} 
\eps^2 {\tA}_{(\pm2)1} (\eps
(x-c_gt),\eps^2t)\,\EE^{\pm2}\\
\eps^2 {\tA}_{(\pm2)2}( \eps
(x-c_gt),\eps^2t)\,\EE^{\pm2}
 \end{array} \right),
\end{align*}
where $0 < \eps \ll 1$, $\EE = e^{i(k_0x-  \omega_0 t)} $, $ \omega_0 = \omega(k_0) $,
${\tA}_{-j} = \overline{\tA}_j$, and 
${\tA}_{-j\ell} = \overline{\tA}_{j\ell}$.   
\begin{remark}{ Our ansatz leads to waves moving to the right.
For waves moving to the left one has to replace in the above ansatz
the vector $
(1 , 0 )^T $ by $
(0 , 1 )^T $ as well as $ -\omega_0$ by $ \omega_0$
and $c_g$ by $-c_g$.
}
\end{remark}
 We insert our ansatz \eqref{ans1c} into system \eqref{diag1}-\eqref{diag2}. Then we
replace the dispersion relation $\omega = \omega(k) $ in all terms of the form $\omega {\tA}_{j}\EE^{j}$ or $\omega {\tA}_{j\ell}\EE^{j}$ by their Taylor expansions 
 around $k=jk_0$. (Details of these expansions are contained in Lemma 25 of \cite{SW10}, for example.) After that, we
equate the coefficients of the $\eps^m\EE^j$ to zero. 

We find that
the coefficients of $\eps \EE^1$ and $\eps^2 \EE^1$ vanish identically
due to the definition of $ \omega$ and 
$ c_g$.
For $\eps^3\EE^1 $ we obtain
$$
\partial_T {\tA}_1 = \frac{1}{2}i\, \partial_k^2 \omega(k_0)\, \partial^2_X {\tA}_1 +
 g_1\,,
$$
where $g_1$ is a sum of multiples of 
$ {\tA}_1 {\tA}_{0\ell} $ and $ {\tA}_{-1} {\tA}_{2\ell} $.
In the next steps we obtain algebraic relations
such that  the $ {\tA}_{2\ell} $ can be expressed in terms of $({\tA}_1)^2$ and 
the $ {\tA}_{0\ell} $ in terms of $ |{\tA}_1|^2 $, respectively.

For $\eps^2 \EE^2$ we obtain
\begin{align*}
( -2 \omega_0 + \omega(2 k_0) ) {\tA}_{21} & =  \gamma_{21} ({\tA}_1)^2\,,
 \\
( -2 \omega_0 -  \omega(2 k_0) ) {\tA}_{22} & =  \gamma_{22} ({\tA}_1)^2\,,
 \nonumber
\end{align*}
with coefficients $ \gamma_{2\ell} \in \C $.
Since   
$-2 \omega_0 \pm \omega(2 k_0) \neq 0 $, which follows from the explicit form of $\omega(k)$, the $ {\tA}_{2\ell} $ are well-defined in terms of $({\tA}_1)^2$.

For $\eps^2\EE^0$ we find 
\begin{align*} 
 \lim\limits_{k \to 0^{-}} \omega(k)\,  {\tA}_{01} & =  \gamma_{01} ({\tA}_1 {\tA}_{-1} )\,,\\ \nonumber
\lim\limits_{k \to 0^{+}} \omega(k)\, {\tA}_{02} & = 
\gamma_{02} ({\tA}_1 {\tA}_{-1} )\,,
\end{align*}
where now $ \gamma_{0\ell} \in \R $ according to the fact that 
we consider a real-valued problem.
Since $\lim\limits_{k \to 0^{\pm}} \omega(k)  \neq 0$, we can express the $ {\tA}_{0\ell} $ in terms of $  |{\tA}_1|^2 $.
 
As mentioned above the nonlinear term in the equation for $\eps^3 \EE^1 $ 
include terms  consisting of combinations of ${\tA}_1$ with 
the ${\tA}_{0\ell}$ and of ${\tA}_{-1}$ with the ${\tA}_{2\ell}$. Eliminating
${\tA}_{0\ell}$ and ${\tA}_{2\ell}$ by the algebraic relations obtained 
for $\eps^2 \EE^0 $
and $\eps^2 \EE^2 $ gives finally the NLS equation
 \begin{equation}\label{nlsderive}
 \partial_T {\tA}_1 = i\, \frac{\omega''(k_0)}{2}\, \partial^2_X {\tA}_1 
+ i\nu_2 (k_0) {\tA}_1 |{\tA}_1|^2\,,
\end{equation}
with a $\nu_2(k_0) \in \R$.

To prove the approximation property of the NLS equation \eqref{nlsderive} it will be helpful to make the residual 
\begin{equation}
{\rm Res}_{u}(\veps \tPsi) =
\left( \ba{c}  {\rm Res}_{u_{-1}}(\veps \tPsi) \\
{\rm Res}_{u_{1}}(\veps \tPsi) \ea \right)
\end{equation}
which contains all terms that do not cancel after inserting ansatz 
\eqref{ans1c} into system \eqref{diag1}-\eqref{diag2}, smaller by modifying 
$\veps \tPsi$ in the following way. First, the above approximation $\veps \tPsi$ is extended by higher order terms. Secondly, by some cut-off function 
the support of the modified approximation 
in Fourier space is restricted to small neighborhoods 
of a finite number of integer multiples of the basic wave number $ k_0 > 0 $.
Since the approximation in Fourier space is strongly concentrated around these 
wave numbers, the approximation is only changed slightly by this modification,
but this second step will give us a simpler control
of the error and makes the approximation an analytic function.

Since $\pm \omega(m k_0) \neq m \omega(k_0)$ 
for all integers $m \geq 2$,
we can proceed analogously as in \cite{DSW12} to
replace $\veps \tPsi$ by a new approximation $\eps\Psi$ of the form
\begin{equation}
\eps\Psi = \eps\Psi_1 + \eps\Psi_{-1} + \eps^2\Psi_q\,,
\end{equation}
where
\begin{align*}
\eps \Psi_{\pm 1} & = \eps \psi_{\pm 1} \!\left( \ba{c} 
1 \\ 0  \ea \right)  =  \eps A_{\pm 1} (\eps
(x-c_g t),\eps^2 t)
\,\EE^{\pm 1} \!\left( \ba{c} 
1 \\ 0  \ea \right),\\[3mm]
\eps^2\Psi_q & = \left(\begin{array}{c}
\eps^2 \psi_{q_{-1}} \\
\eps^2 \psi_{q_{1}}
\\ \end{array}\right) = 
 \eps^2\Psi_0
 + \eps^2\Psi_2  + \eps^2\Psi_{-2}  + \eps^2\Psi_h\,,\\[3mm]
\eps^2 \Psi_0 & =  \left(\begin{array}{c}\eps^2 \psi_{01}\\
\eps^2 \psi_{02} \\
\end{array}\right)  =  \left(\begin{array}{c}\eps^2 A_{01} (\eps
(x-c_gt),\eps^2t)\\
\eps^2 A_{02} (\eps (x-c_gt),\eps^2t)\\
\end{array}\right),
\\[3mm]
\eps^2 \Psi_{\pm2} & = \left(\begin{array}{c}\eps^2 \psi_{(\pm2)1}\\
\eps^2 \psi_{(\pm2)2} 
\\ \end{array}\right)  =  \left(\begin{array}{c} 
\eps^2 A_{(\pm2)1} (\eps (x-c_gt),\eps^2t)\,\EE^{\pm2}\\
\eps^2 A_{(\pm2)2} (\eps
(x-c_gt),\eps^2t)\,\EE^{\pm2}\\
\end{array} \right),\\[3mm]
\eps^2\Psi_h =\, &
\sum\limits_{j=-1,1\atop n=1,2,3} \left(
\begin{array}{c} \eps^{1+n}A^n_{j1} (\eps(x -c_gt),\eps^2t)\EE^j\\
\eps^{1+n}A^{n}_{j2} (\eps(x -c_gt),\eps^2t)\EE^j\\
\end{array}\right)\\
& + \sum\limits_{j=-2,2 \atop n=1,2} 
\left(\begin{array}{c}\eps^{2+n}A^n_{j1} (\eps(x -c_gt),\eps^2t)\EE^j\\
\eps^{2+n}A^n_{j2} (\eps(x -c_gt),\eps^2t)\EE^j\\
\end{array}\right)\\
& + \sum\limits_{n=1,2} 
\left(\begin{array}{c}\eps^{2+n}A^n_{01} (\eps(x -c_gt),\eps^2t)\\
\eps^{2+n}A^n_{02} (\eps(x -c_gt),\eps^2t)
\end{array}\right) \\
\,& + \sum\limits_{j=-3,3 \atop n=0,1} 
\left(\begin{array}{c}\eps^{3+n}A^n_{j1} (\eps(x -c_gt),\eps^2t)\EE^j\\
\eps^{3+n}A^n_{j2} (\eps(x -c_gt),\eps^2t)\EE^j\\
\end{array}\right)\, +
\end{align*}
\begin{align*}
\,& + \sum\limits_{j=-4,4} 
\left(\begin{array}{c}\eps^{4} A^n_{j1} (\eps(x -c_gt),\eps^2t)\EE^j\\
\eps^{4} A^n_{j2} (\eps(x -c_gt),\eps^2t)\EE^j\\
\end{array}\right),
\end{align*}
${A}_{-j} = \overline{A}_j$, and 
${A}_{-j\ell} = \overline{A}_{j\ell}$, 
which has compact support in Fourier space for all $0 < \eps \ll 1$.
Then, exactly as in Section 2 of \cite{DSW12}, 
the following estimates for the modified residual hold.
\begin{lemma} \label{lem2}
Let $ s_A \geq {6} $ and
$ {\tA}_1 \in C([0,T_0], H^{s_A}(\R,\C)) $ be a solution
of the NLS equation \eqref{nlsderive}  with 
$$ \sup_{T \in [0,T_0]} \| 
\tA_{1} \|_{H^{s_A}} \leq C_A . $$ 
Then for all $s \geq 0$ there exist $ C_{Res}, C_{\Psi}, \eps_0>0  $ depending on  $C_A$ such that for all
$ \eps \in (0,\eps_0) $ the approximation $\eps \Psi$ satisfies
\begin{align} \label{RES1}
\sup_{t \in [0,T_0/\eps^2]} \|  {\rm Res}_u(\eps \Psi)
\|_{H^s}
& \leq  C_{Res}\, \eps^{{9/2}}, \\ \label{RES2}
\sup_{t \in [0,T_0/\eps^2]} \|\eps \Psi - (\eps {\tPsi_{1}} + \eps {\tPsi_{-1}})
\|_{H^{{s_A}}}
& \leq  {C_{\Psi}}\, \eps^{3/2},
\\
\label{RES3}
\sup_{t \in [0,T_0/\eps^2]} (\|\widehat{\Psi}_{\pm 1} \|_{L^1({s+1})(\R,\C)}
+ \|\widehat{\Psi}_q \|_{L^1({s+1})(\R,\C)})
& \leq {C_{\Psi}}\,. 
\end{align}
\end{lemma}

The proof of Lemma \ref{lem2} goes analogously as the proof of Lemma 2.6 in \cite{DSW12}. The fact that for all $s \geq 0$ the first and the third estimate are valid for appropriate constants $ C_{Res}$ and $C_{\Psi}$ is a consequence 
of the fact that our approximation $\eps \Psi$ has compact support in Fourier 
space. The approximation $\eps \Psi$ differs so slightly from the actual NLS approximation $\eps ({\tPsi_{1}} +  {\tPsi_{-1}})$ and higher order asymptotic expansions of the exact solution, which are
needed to make the residual sufficiently small, that the bounds in \eqref{RES1}-\eqref{RES2} hold
if $s_A \geq 6$. This is shown with the help of the estimate
\begin{equation*} 
 \| (\chi_{[-\delta,\delta]}-1)\, \eps^{-1} 
\widehat{f} ( \eps^{-1} \cdot) \|_{L^2(m)}  \lesssim 
\eps^{m+{M}-1/2} \| f \|_{H^{m+{M}}}
\end{equation*}
for all $M,m \geq 0$, where $\chi_{[-\delta,\delta]}$ is the characteristic function on $[-\delta,\delta]$. 

\begin{remark} \label{remneucam2}
{\rm 
The bound \eqref{RES3} will be  used for instance to estimate 
$$ \| \psi_{j}f \|_{H^{s}} \leq C \| \psi_{j}  \|_{C^{s}_b}\| f 
\|_{H^{s}} \leq  C \| \widehat{\psi}_{j}  \|_{L^1(s)(\R,\C)}\| f
\|_{H^{s}} $$
without loss of powers in $ \eps $ as it would be the case with $\| \widehat{\psi}_{j}  \|_{L^2(s)(\R,\C)}$.}
\end{remark}

Moreover, by an analogous argumentation as in the proof of Lemma 3.3 in \cite{DSW12} we obtain the fact that $\partial_t \psi_{\pm 1}$ can be approximated by
$-i \omega \psi_{\pm 1}$. More precisely, we get 
\begin{lemma}\label{deriv_psi} For all $s>0$
there exists a constant $C_{\psi} > 0$ such that
\begin{equation} \label{dtpsi}
\| \partial_t \widehat{\psi}_{\pm 1} + i \omega \widehat{\psi}_{\pm 1} \|_{L^{1}(s)} \le C_{\psi}\, \eps^2\,.
\end{equation}
\end{lemma}

\section{The Error Estimates}
\label{sec:err}
Now, we write $u_{\pm1}$ as approximation plus error:
\begin {equation}
 \left(
\begin{array}{c} \widehat{u}_{-1} \\ \widehat{u}_{1}
\end{array}
\right)  = \veps \Psi + \veps^{5/2}  \left(
\begin{array}{c} \widehat{R}_{-1} \\ \widehat{R}_{1}
\end{array}
\right)\,.
\end{equation}
This yields
\begin{align} \label{err-1}
\partial_t \widehat{R}_{-1} (k,t)  = &   -i \omega(k) \widehat{R}_{-1}(k,t)
 - \veps i \rho(k) (\widehat{\psi} \ast (\widehat{R}_{-1}+\widehat{R}_1))(k,t)\\ & - \frac{1}{2} \veps^{5/2} i \rho(k) (\widehat{R}_{-1}+\widehat{R}_1)^{\ast 2}(k,t) + \veps^{-5/2} \widehat{{\rm Res}}_{u_{-1}}(\veps \Psi)(k,t)\,, \nonumber\\[2mm]
\partial_t \widehat{R}_1(k,t)  = & \;  i \omega(k) \widehat{R}_1(k,t) + \veps i \rho(k) (\widehat{\psi} \ast (\widehat{R}_{-1} + \widehat{R}_1))(k,t) \label{err1}\\ &
+ \frac{1}{2} \veps^{5/2} i \rho(k) (\widehat{R}_{-1}+\widehat{R}_1)^{\ast 2}(k,t) + \veps^{-5/2} \widehat{{\rm Res}}_{u_{1}}(\veps \Psi)(k,t) \,, \nonumber
\end{align}
where $\psi=\psi_{-1}+\psi_{1}+ \eps \psi_{q_{-1}} + \eps \psi_{q_{1}}$.
\medskip

In order to control the error we use the energy 
\begin{equation} 
\mathcal{E}_s = \sum_{\ell=0}^s E_{\ell}\,, 
\end{equation} 
\begin{equation*}
E_{\ell} = \sum_{j_1\in \{\pm1\}} \Big( \frac{1}{2}   \int_{\R} (\partial_x^{\ell} R_{j_1})^2\,dx  +     \veps \sum_{j_2 \in \{\pm1\}} \int_{\R}\partial_x^{\ell} R_{j_1} \partial_x^{\ell} N_{j_1 j_2}(\psi,R_{j_2})\,dx \Big), 
\end{equation*}
with
\begin{equation*}
\widehat{N}_{j_1 j_2}(\psi,R_{j_2})(k) =  \int_{\R} \widehat{n}_{j_1 j_2}(k,k-m,m) \widehat{\psi}(k-m) \widehat{R}_{j_2}(m)\,dm\,,  
\end{equation*}
\begin{equation*}
\widehat{n}_{j_1 j_2}(k,k-m,m) = \frac{-j_1 \rho(k)\, \chi(k-m)}{-j_1\omega(k)-\omega(k-m)+j_2\omega(m)}\,,
\end{equation*}
where $s=s_A \geq 6$ and $\chi$ is the characteristic function on $\rm{supp}\,\widehat{\psi}$.

\begin{lemma} \label{lem31}
The operators $N_{j_1 j_2}$ have the following properties:\\[2mm]
{\bf a)}\; Fix $h\in L^2(\R,\R)$. Then $f \mapsto N_{jj}(h,f)$ defines a continuous linear map from $H^1(\R,\R)$ into $L^2(\R,\R)$ and $f \mapsto N_{j-j}(h,f)$ a continuous linear map from $L^2(\R,\R)$ into $L^2(\R,\R)$. In particular,
for all $f\in H^1(\R,\R)$ we have
\begin{align} \label{as1}
N_{jj}(h,f)&= -j\partial_x(G_{jj}h\,f) + Q_{jj}(h,f)\,,\\[2mm] \label{as2}
N_{j-j}(h,f)&= G_{j-j}h\,f + Q_{j-j}(h,f)\,,
\end{align}
with
\begin{align*}
\widehat{G_{jj}h}(k)&=  \frac{\chi(k)}{-i(\omega(k)+jk)}\widehat{h}(k)\,,\\[2mm]
\widehat{G_{j-j}h}(k)&= \frac{1}{2} \chi(k) \widehat{h}(k)\,,\\[2mm]
\|Q_{j\pm j}(h,f)\|_{H^1}&= \mathcal{O}(\|h\|_{L^2} \|f\|_{L^2})\,.
\end{align*}
{\bf b)}\; For all $f \in H^1(\R,\R)$ we have
\begin{align} \label{nf}
-j_1 i \omega N_{j_1 j_2}(\psi,f)- N_{j_1 j_2}(i \omega \psi,f)+ j_2 N_{j_1 j_2}(\psi,i \omega f) &= -j_1 i \rho (\psi f)\,,
\end{align}
where the operators $\omega$ and $\rho$ are defined by the symbols \eqref{omegasym}-\eqref{rhosym}.\\[2mm]
{\bf c)}\;
For all $f,g,h \in H^1(\R,\R)$ we have
\begin{align} \label{partN}
\int_{\R} f\, N_{j_1j_2}(h,g)\,dx &= -\frac{j_1}{j_2} \int_{\R} N_{j_2j_1}(h,f)\, g\,dx + \int_{\R} S_{j_2j_1}(\partial_x h,f)\, g\,dx\,,
\end{align}
where
\begin{align*}
\widehat{S}_{j_2j_1}(\partial_x h,f)(k)&=  \int_{\R} \widehat{{s}}_{j_2j_1}(k,k-m,m) \widehat{\partial_x h}(k-m) \widehat{f}(m)\,dm\,,  
\end{align*}
with
\begin{align*}
\widehat{{s}}_{j_2j_1}(k,k-m,m) &= \frac{-j_1\, (\rho(k)-\rho(m))\, \chi(k-m)}{(k-m)\,i\,(-j_2 \omega(k)-\omega(k-m)+j_1 \omega(m))}\,.
\end{align*}
In particular, we have
\begin{align}  \label{as3}
S_{jj}(\partial_x h,f)&= -jG_{jj}\partial_x h\,f + \tilde{Q}_{jj}(\partial_x h,f)\,,
\end{align}
with
\begin{align*}
\|\tilde{Q}_{jj}(\partial_x h,f)\|_{H^2}&= \mathcal{O}(\|h\|_{L^2} \|f\|_{L^2})\,.
\end{align*}
\end{lemma}
{\bf Proof.}
Since $\rm{supp}\,\widehat{\psi}$ is compact by construction, there exists a $k_1 >0$ with 
${\rm {supp}}\,\chi \subset [-k_1,k_1]$. Therefore, the non-resonance condition
\begin{equation}
\inf_{{k\in \R,\atop p \in [-k_1,k_1],}  \atop  j_1, j_2 \in \{\pm1\}}\, | -j_1\omega(k)-\omega(p)+j_2\omega(k-p)| \geq  C > 0
\end{equation}
is satisfied with a constant $C=C(k_1) \rightarrow 0$ for $k_1 \rightarrow \infty$, see \cite{Schn05}, which implies $|\widehat{n}_{j_1 j_2}(k,k-m,m)| < \infty$ for all $k,m \in \R$.

Next, we analyze the asymptotic behavior of the $\widehat{n}_{j_1 j_2}(k,k-m,m)$ for $|k| \rightarrow \infty$. 
We have
\begin{align} \label{asomega}
\omega(k) &= k + \mathcal{O}(|k|^{-1}) \qquad \mathrm{for}\;  |k| \rightarrow \infty\,,\\[2mm] \label{asdkomega}
\omega'(k) &= 1 + \mathcal{O}(k^{-2}) \qquad \mathrm{for}\;  |k| \rightarrow \infty\,,\\[2mm]
\label{asrho}
\rho(k) &= k + \mathcal{O}(|k|^{-1}) \qquad \mathrm{for}\;  |k| \rightarrow \infty\,.
\end{align}
By the mean value theorem we get
\begin{align*}
\widehat{n}_{j j}(k,k-m,m) &=  \frac{j\rho(k)\, \chi(k-m)}{\omega(k-m)+j(\omega(k)-\omega(m))}\\[2mm]
&= \frac{j\rho(k)\, \chi(k-m)}{\omega(k-m)+j(k-m)\,\omega'(k-\vartheta (k,m)(k-m))}\,, 
\end{align*}
with $\vartheta(k,m) \in [0,1]$.
Using again the fact that $\rm{supp}\,\chi$ is compact, we conclude with the help of the expansions \eqref{asdkomega} and \eqref{asrho} that
\begin{align*}
\widehat{n}_{j j}(k,k-m,m) 
&= \frac{j\, (k + \mathcal{O}(|k|^{-1}))\,\chi(k-m)}{\omega(k-m)+j(k-m)(1+ \mathcal{O}(k^{-2}))} \qquad \mathrm{for}\;  |k| \rightarrow \infty\,,\\[2mm]
&= \Big(\frac{jk}{\omega(k-m)+j(k-m)} + \mathcal{O}(|k|^{-1})\Big)\,\chi(k-m)\qquad \mathrm{for}\;  |k| \rightarrow \infty\,.
\end{align*}
Exploiting once more the compactness of  $\rm{supp}\,\chi$ as well as the expansions \eqref{asomega}-\eqref{asrho} yields
\begin{align*}
\widehat{n}_{j -j}(k,k-m,m) &= \frac{j\rho(k)\, \chi(k-m)}{j(\omega(k)+\omega(k-(k-m)))+ \omega(k-m)}\\[2mm]
&= \frac{\rho(k)\, \chi(k-m)}{2\omega(k)(1+ \mathcal{O}(|k|^{-1}))} \qquad \mathrm{for}\;  |k| \rightarrow \infty\,,\\[2mm]
&= \Big(\frac{1}{2} + \mathcal{O}(|k|^{-1})\Big)\,\chi(k-m) \qquad \mathrm{for}\;  |k| \rightarrow \infty\,.
\end{align*}
These asymptotic expansions of the $\widehat{n}_{j_1 j_2}(k,k-m,m)$ imply \eqref{as1}-\eqref{as2}. 

Finally, since
\begin{equation*}
\widehat{n}_{j_1 j_2}(-k,-(k-m),-m)= \widehat{n}_{j_1j_2}(k,k-m,m) \in \R
\end{equation*}
and $\psi$ is real-valued, we obtain the validity of all assertions of a).
\medskip

b) is a direct consequence of the construction of the operators $N_{j_1j_2}$.
\medskip

In order to prove c) we compute for all $f,g,h \in H^1(\R,\R)$:
\begin{align*}
&
\int_{\R} \overline{\widehat{f}(k)}\, \widehat{N}_{j_1j_2}(h,g)(k)\,dk \\[2mm]
=\;& \int_{\R} \int_{\R}  \overline{\widehat{f}(k)}\;  \frac{-j_1 \rho(k)\, \chi(k-m)}{-j_1 \omega(k)-\omega(k-m)+j_2 \omega(m)}\; \widehat{h}(k-m)\, \widehat{g}(m)\,dm\,dk\\[3mm]
=\;& \int_{\R} \int_{\R}  \overline{\widehat{g}(-m)}\;  \frac{-j_1 \rho(k)\, \chi(k-m)}{j_2 \omega(m)-\omega(k-m)-j_1 \omega(k)}\; \widehat{h}(k-m)\, \widehat{f}(-k)\,dk\,dm\\[3mm]
=\;& \int_{\R} \int_{\R}  \overline{\widehat{g}(k)}\;  \frac{j_1 \rho(m)\, \chi(k-m)}{-j_2 \omega(k)-\omega(k-m)+j_1 \omega(m)}\; \widehat{h}(k-m)\, \widehat{f}(m)\,dm\,dk\\[3mm]
=\;& \int_{\R} \int_{\R}  \overline{\widehat{g}(k)} \; \frac{j_1 \rho(k)\, \chi(k-m)}{-j_2 \omega(k)-\omega(k-m)+j_1 \omega(m)}\; \widehat{h}(k-m)\, \widehat{f}(m)\,dm\,dk\\[3mm]
&+ \int_{\R} \int_{\R}  \overline{\widehat{g}(k)}\;  \frac{-j_1 (\rho(k)-\rho(m))\, \chi(k-m)}{-j_2 \omega(k)-\omega(k-m)+j_1 \omega(m)}\; \widehat{h}(k-m)\, \widehat{f}(m)\,dm\,dk\\[3mm]
=\;& -\frac{j_1}{j_2} \int_{\R} \overline{\widehat{g}(k)}\, \widehat{N}_{j_2j_1}(h,f)(k)\,dk
+ \int_{\R} \overline{\widehat{g}(k)}\, \widehat{S}_{j_2j_1}(\partial_x h,f)(k)\,dk \,,
\end{align*}
which yields \eqref{partN}, and due to \eqref{as1} and \eqref{asrho} we obtain \eqref{as3}.
\qed
\medskip

The assertions of Lemma \ref{lem31} a), c) and the Cauchy-Schwarz inequality imply
\begin{corollary} \label{cor32}
$\sqrt{\mathcal{E}_s}$ is equivalent to $\|R_1\|_{H^s}+\|R_{-1}\|_{H^s}$ for sufficiently small $\veps>0$.
\end{corollary}

Since the right-hand sides of the error equations \eqref{err-1}-\eqref{err1} lose one derivative, we will need the following identities to control the time evolution of $\mathcal{E}_s$. 

\begin{lemma} \label{lem33}
Let $j \in \{\pm1\}$, $a_j \in H^2(\R,\R)$, and $f_j \in H^1(\R,\R)$. Then we have
\begin{align}
\label{dx1} 
\int_{\R} a_{j}\, f_{j}\, \partial_x f_{j}\, dx
=& -\frac{1}{2}\, \int_{\R} \partial_xa_{j}\, f_{j}^2\, dx\,,\\[3mm]
\label{dx2}
\sum_{j \in \{\pm1\}}\, \int_{\R} a_{j}\, f_{j}\, \partial_x f_{-j}\,dx
=&\; \frac{1}{2}\, \int_{\R} (a_{-1}-a_1)\,(f_1+f_{-1})\, \partial_x(f_1-f_{-1})\,dx \\[2mm]
\nonumber
& + \mathcal{O} \big((\|a_1\|_{H^2}+ \|a_{-1}\|_{H^2})(\|f_1\|^2_{L^2}+\|f_{-1}\|^2_{L^2})\big)\,.
\end{align}
\end{lemma}
{\bf Proof.}
Identity \eqref{dx1} follows directly by partial integration. 
Using again partial integration, the Cauchy-Schwarz inequality, and \eqref{dx1} we obtain 
\begin{align*}
& \sum_{j \in \{\pm1\}}\, \int_{\R} a_{j}\, f_{j}\, \partial_x f_{-j}\,dx\\[2mm]
 =\;& \frac{1}{2}\, \sum_{j \in \{\pm1\}} \Big(\int_{\R} a_{j}\, f_{j}\, \partial_x f_{-j}\,dx - \int_{\R} a_{j}\, \partial_x f_{j}\, f_{-j}\,dx - \int_{\R} \partial_xa_{j}\, f_{j}\, f_{-j}\,dx \Big) \\[2mm] 
 =\;& \frac{1}{2}\, \Big( \int_{\R} (a_{-1}-a_1)\, f_{-1}\, \partial_x f_{1}\,dx - \int_{\R} (a_{-1}-a_1)\, f_1\, \partial_x f_{-1}\,dx \Big)\\[2mm]  
& + \mathcal{O} \big((\|a_1\|_{H^2}+ \|a_{-1}\|_{H^2})(\|f_1\|^2_{L^2}+\|f_{-1}\|^2_{L^2})\big)\\[2mm]
 =\;& \frac{1}{2}\, \Big( \int_{\R} (a_{-1}-a_1)\, (f_1+f_{-1})\, \partial_x f_{1}\,dx - \int_{\R} (a_{-1}-a_1)\,  (f_1+f_{-1})\, \partial_x f_{-1}\,dx\Big)\\[2mm]  
& + \mathcal{O} \big((\|a_1\|_{H^2}+ \|a_{-1}\|_{H^2})(\|f_1\|^2_{L^2}+\|f_{-1}\|^2_{L^2})\big)\\[2mm]
=\;& \frac{1}{2}\, \int_{\R} (a_{-1}-a_1)\,(f_1+f_{-1})\, \partial_x(f_1-f_{-1})\,dx \\[2mm]
& + \mathcal{O} \big((\|a_1\|_{H^2}+ \|a_{-1}\|_{H^2})(\|f_1\|^2_{L^2}+\|f_{-1}\|^2_{L^2})\big)\,.
\end{align*}
\qed
\medskip

Now, we are prepared to analyze $\partial_t E_{\ell}$. We compute
\begin{align*}
\partial_t E_{\ell} =& \sum_{j_1 \in \{\pm1\}} \Big( \int_{\R} \partial^{\ell}_x R_{j_1}\, \partial_t\partial^{\ell}_x R_{j_1} \,dx\,
+   \veps \sum_{j_2 \in \{\pm1\}} \Big( \int_{\R} \partial_t\partial^{\ell}_x R_{j_1}\, \partial^{\ell}_x N_{j_1j_2}(\psi,R_{j_2})\,dx
\\[2mm] & 
+ \int_{\R} \partial^{\ell}_x R_{j_1}\, \partial^{\ell}_x N_{j_1j_2}(\psi,\partial_t R_{j_2})\,dx\,
+ \int_{\R} \partial^{\ell}_x R_{j_1}\, \partial^{\ell}_x N_{j_1j_2}(\partial_t \psi, R_{j_2} )\,dx \Big)\Big)\,.
\end{align*}
Using the error equations \eqref{err-1}-\eqref{err1} we get
\begin{align*}
\partial_t E_{\ell} 
=\;& \sum_{j_1 \in \{\pm1\}} \Big(\, j_1 \int_{\R} \partial^{\ell}_x R_{j_1}\, i\omega \partial^{\ell}_x R_{j_1} \,dx\\[1mm] &
\qquad\qquad\! + \int_{\R} \partial^{\ell}_x R_{j_1}\, \veps^{-5/2} \partial^{\ell}_x \mathrm{Res}_{u_{j_1}}\!(\veps\Psi)\,dx\, \Big)
\\[2mm] 
& 
+\; \veps \sum_{j_1, j_2 \in \{\pm1\}}\! \Big(\, j_1 \int_{\R} \partial^{\ell}_x R_{j_1}\, i \rho \partial^{\ell}_x({\psi} R_{j_2})\,dx \\[1mm]
& 
\qquad\qquad\qquad\,\,\, + j_1  \int_{\R} i \omega \partial^{\ell}_x R_{j_1}\, \partial^{\ell}_xN_{j_1j_2}(\psi, R_{j_2})\,dx
\\[2mm]
&\qquad\qquad\qquad\,\,\, + j_2  \int_{\R} \partial^{\ell}_x R_{j_1}\, \partial^{\ell}_x N_{j_1j_2}(\psi, i \omega R_{j_2})\,dx
\\[2mm]
&\qquad\qquad\qquad\,\,\, - \int_{\R} \partial^{\ell}_x R_{j_1}\, \partial^{\ell}_x N_{j_1j_2}(i \omega \psi, R_{j_2})\,dx
\\[2mm]
&\qquad\qquad\qquad\,\,\,
+ \int_{\R} \partial^{\ell}_x R_{j_1}\, \partial^{\ell}_x N_{j_1j_2}(\partial_t \psi +i \omega \psi, R_{j_2})\,dx\\[2mm] 
& \qquad\qquad\qquad\,\,\,
- \int_{\R} \veps^{-5/2} \partial^{\ell+1}_x \mathrm{Res}_{u_{j_1}}\!(\veps\Psi)\, \partial^{\ell-1}_x N_{j_1j_2}(\psi,R_{j_2})\,dx
\\[2mm]& \qquad\qquad\qquad\,\,\,
+ \int_{\R} \partial^{\ell}_x R_{j_1}\, \partial^{\ell}_x N_{j_1j_2}(\psi,\veps^{-5/2} \mathrm{Res}_{u_{j_2}}\!(\veps\Psi))\,dx\, \Big)\\[2mm]
&
+\; \veps^{2}  \sum_{j_1,j_2,j_3 \in \{\pm1\}}\!  \Big(\, j_1\int_{\R}  i \rho \partial^{\ell}_x(\tilde{\psi} R_{j_3})\, \partial^{\ell}_x N_{j_1j_2}(\psi, R_{j_2})\,dx\\[1mm]
&\qquad\qquad\qquad\qquad\! +j_2 \int_{\R} \partial^{\ell}_x R_{j_1}\, \partial^{\ell}_x N_{j_1j_2}(\psi, i\rho (\tilde{\psi} R_{j_3}))\,dx \,\Big) \\[2mm]&
+\;\veps^{5/2} \sum_{j_1, j_2,j_3  \in \{\pm1\}} \frac{j_1}{2}  \int_{\R} \partial^{\ell}_x R_{j_1}\, i \rho \partial_x^{\ell}(R_{j_3} R_{j_2})\,dx\,,
\end{align*}
where $\tilde{\psi}=\psi + \frac{1}{2}\veps^{3/2} (R_{1}+R_{-1})$.

Due to the skew symmetry of $i\omega$ the first integral equals zero. Since the operators $N_{j_1j_2}$ satisfy \eqref{nf}, the third integral cancels with the sum of
the fourth, the fifth, and the sixth integral. 

Moreover, because of the estimates \eqref{RES1} and \eqref{RES3} for the residual, the bound \eqref{dtpsi} for $\partial_t \psi_{\pm 1} + i \omega \psi_{\pm 1}$, the regularity properties of the operators $N_{j_1j_2}$ from Lemma \ref{lem31} a), identity \eqref{dx1}, and Corollary \ref{cor32}, the second, the seventh, the eighth, and the ninth integral can be bounded by $C\veps^{2}(\mathcal{E}_s+1)$ for a constant $C>0$. Hence, we have
\begin{align*}
\partial_t E_{\ell} =\;& \veps^{2} \sum_{j_1,j_3  \in \{\pm1\}} \Big(\,j_1 \int_{\R}  i \rho \partial^{\ell}_x(\tilde{\psi} R_{j_3})\, \partial_x^{\ell}N_{j_1j_1}(\psi, R_{j_1})\,dx\\[1mm]
&\qquad\qquad\qquad+ j_1 \int_{\R} \partial^{\ell}_x R_{j_1}\, \partial^{\ell}_xN_{j_1j_1}(\psi, i\rho (\tilde{\psi} R_{j_3}))\,dx 
\\[2mm]&
\qquad\qquad\qquad+j_1 \int_{\R}  i \rho \partial^{\ell}_x(\tilde{\psi} R_{j_3})\, \partial_x^{\ell}N_{j_1-j_1}(\psi, R_{-j_1})\,dx\; \\[2mm]
&\qquad\qquad\qquad- j_1 \int_{\R} \partial^{\ell}_x R_{j_1}\, \partial^{\ell}_xN_{j_1-j_1}(\psi, i\rho (\tilde{\psi} R_{j_3}))\,dx\, \Big) \\[2mm]
& +\,\veps^{5/2} \sum_{j_1, j_2,j_3  \in \{\pm1\}} \frac{j_1}{2}   \int_{\R} \partial^{\ell}_x R_{j_1}\, i \rho \partial_x^{\ell}(R_{j_3} R_{j_2})\,dx \\[2mm] 
& +\, \veps^{2}\, \mathcal{O}(\mathcal{E}_s+1)\\[2mm] 
=:\;& \sum_{j=1}^5 I_j + \veps^{2}\,\mathcal{O}(\mathcal{E}_s+1)\,.
\end{align*}
First, we analyze $I_1+I_2$. To extract all terms with more than $\ell$ spatial derivatives falling on $R_1$ or $R_{-1}$ we apply Leibniz's rule and get
\begin{align*}
I_1+ I_2 
=\;& \veps^{2} \sum_{j_1, j_3 \in \{\pm1\}} \Big(\, j_1 \int_{\R}  i \rho \partial^{\ell}_x (\tilde{\psi} R_{j_3})\, N_{j_1j_1}(\psi, \partial^{\ell}_x R_{j_1})\,dx\\[1mm]
&\qquad\qquad\qquad+ \ell j_1 \int_{\R}  i \rho \partial^{\ell}_x(\tilde{\psi} R_{j_3})\, N_{j_1j_1}(\partial_x \psi, \partial^{\ell-1}_x R_{j_1})\,dx\\[2mm]
&\qquad\qquad\qquad+j_1 \int_{\R} \partial^{\ell}_x R_{j_1} \,N_{j_1j_1}(\psi, i \rho \partial^{\ell}_x(\tilde{\psi} R_{j_3}))\,dx \\[2mm]
&\qquad\qquad\qquad+\ell j_1 \int_{\R} \partial^{\ell}_x R_{j_1}\, N_{j_1j_1}(\partial_x\psi, i \rho\partial^{\ell-1}_x (\tilde{\psi} R_{j_3}))\,dx \,\Big)\\[2mm]
&+\, \veps^{2}\,\mathcal{O}(\mathcal{E}_s+ \veps^{3/2} \mathcal{E}^{3/2}_s)\,. 
\end{align*}
Because of \eqref{partN} we obtain
\begin{align*}
I_1+ I_2  =\;& \veps^{2} \sum_{j_1, j_3 \in \{\pm1\}} \Big(\, j_1 \int_{\R}  i \rho \partial^{\ell}_x(\tilde{\psi} R_{j_3})\, {S}_{j_1j_1}(\partial_x \psi, \partial^{\ell}_x R_{j_1})\,dx\\[1mm]
&\qquad\qquad\qquad +2\ell j_1 \int_{\R}  i \rho \partial^{\ell}_x(\tilde{\psi} R_{j_3})\, N_{j_1j_1}(\partial_x \psi, \partial^{\ell-1}_x R_{j_1})\,dx\, \Big)\\[2mm]
&+\, \veps^{2}\,\mathcal{O}(\mathcal{E}_s+ \veps^{3/2} \mathcal{E}^{3/2}_s)\,. 
\end{align*}
Using the asymptotic expansions \eqref{as1}, \eqref{as3}, and \eqref{asrho} yields
\begin{align*}
I_1+ I_2 
=\,& -(2\ell+1)\, \veps^{2} \sum_{j_1, j_3 \in \{\pm1\}} \int_{\R} G_{j_1j_1}\!\partial_x \psi\, \tilde{\tilde{\psi}}\, \partial^{\ell}_x R_{j_1}\, \partial^{\ell+1}_x R_{j_3}\,dx \\[2mm]
&+\, \veps^{2}\,\mathcal{O}(\mathcal{E}_s+\veps^{3/2} \mathcal{E}^{3/2}_s)\,,
\end{align*}
where $\tilde{\tilde{\psi}}=\psi + \veps^{3/2} (R_{1}+R_{-1})$.
With the help of \eqref{dx1}, \eqref{dx2}, and
\begin{equation*}
(\widehat{G}_{-1-1}-\widehat{G}_{11})(k) = \frac{2ik\,\chi(k)}{\omega^{2}(k)-k^{2}} = 2ik\, \chi(k)
\end{equation*}
we obtain
\begin{align*}
I_1+ I_2 
=\,& -\frac{2\ell+1}{2}\, \veps^{2}\, \int_{\R} (G_{-1-1}-G_{11}) \partial_x \psi\, \tilde{\tilde{\psi}}\, \,  \partial^{\ell}_x (R_{1}+R_{-1})\, \partial^{\ell+1}_x (R_{1}-R_{-1})\,dx \\[2mm]
&+\, \veps^{2}\,\mathcal{O}(\mathcal{E}_s+\veps^{3/2} \mathcal{E}^{3/2}_s)\\[2mm] 
=\,& -(2\ell+1)\, \veps^{2}\, \int_{\R}\partial^{2}_x \psi\, \tilde{\tilde{\psi}}\, \partial^{\ell}_x (R_{1}+R_{-1})\, \partial^{\ell+1}_x (R_{1}-R_{-1})\,dx \\[2mm]
&+\, \veps^{2}\,\mathcal{O}(\mathcal{E}_s+\veps^{3/2} \mathcal{E}^{3/2}_s)\,,
\end{align*}
and because of
\begin{equation} \label{dxdt}
\partial_t (R_{1}+R_{-1}) = i\omega (R_{1}-R_{-1}) + \veps^{-5/2} ({\rm Res}_{u_{1}}(\veps \Psi) + {\rm Res}_{u_{-1}}(\veps \Psi))
\end{equation}
and \eqref{asomega} we arrive at
\begin{align*}
I_1+ I_2  =\,& -(2\ell+1)\, \veps^{2}\, \int_{\R} \partial^{2}_x \psi\, \tilde{\tilde{\psi}}\, \partial^{\ell}_x (R_{1}+R_{-1})\, \partial_t\partial^{\ell}_x (R_{1}+R_{-1})\,dx \\[2mm]
&+\, \veps^{2}\,\mathcal{O}(\mathcal{E}_s+\veps^{3/2} \mathcal{E}^{3/2}_s)\\[2mm]
=\,& -\frac{2\ell+1}{2}\,\veps^{2}\, \partial_t \int_{\R} \partial^{2}_x \psi\, (\psi + \veps^{3/2} (R_{1}+R_{-1}))\, (\partial^{\ell}_x (R_{1}+R_{-1}))^{2}\,dx \\[2mm]
&+\, \veps^{2}\,\mathcal{O}(\mathcal{E}_s+\veps^{3/2} \mathcal{E}^{3/2}_s)\,.
\end{align*}
The terms $I_3$, $I_4$ and $I_5$ can be analyzed in the same way. Using again Leibniz's rule and \eqref{partN} as well as the asymptotic expansions \eqref{as2}, \eqref{asomega}, and \eqref{asrho} to extract in $I_3+I_4$ all integral terms containing factors with more than $\ell$ spatial derivatives falling on $R_1$ or $R_{-1}$
we get
\begin{align*}
I_3+ I_4 
=\;& \veps^{2} \sum_{j_1, j_3 \in \{\pm1\}} \Big(\, j_1 \int_{\R}  i \rho \partial^{\ell}_x (\tilde{\psi} R_{j_3})\, N_{j_1-j_1}(\psi, \partial^{\ell}_x R_{-j_1})\,dx\\[1mm]
&\qquad\qquad\qquad-j_1 \int_{\R} \partial^{\ell}_x R_{j_1} \,N_{j_1-j_1}(\psi, i \rho \partial^{\ell}_x(\tilde{\psi} R_{j_3}))\,dx\,\Big) \\[2mm]
&+\, \veps^{2}\,\mathcal{O}(\mathcal{E}_s+\veps^{3/2} \mathcal{E}^{3/2}_s)\\[2mm] 
=\;& \veps^{2} \sum_{j_1, j_3 \in \{\pm1\}} \Big(\, j_1 \int_{\R}  i \rho \partial^{\ell}_x (\tilde{\psi} R_{j_3})\, N_{j_1-j_1}(\psi, \partial^{\ell}_x R_{-j_1})\,dx\\[1mm]
&\qquad\qquad\qquad-j_1 \int_{\R} i \rho \partial^{\ell}_x (\tilde{\psi} R_{j_3})\, N_{-j_1j_1}(\psi, \partial^{\ell}_x R_{j_1})\,dx\\[2mm]
&\qquad\qquad\qquad +j_1 \int_{\R}  \partial^{\ell}_x (\tilde{\psi} R_{j_3})\, i\rho {S}_{-j_1j_1}(\partial_x\psi, \partial^{\ell}_x R_{j_1})\,dx\,\Big)\\[2mm]
&+\, \veps^{2}\,\mathcal{O}(\mathcal{E}_s+\veps^{3/2} \mathcal{E}^{3/2}_s)\\[2mm]
=\;& -\veps^{2} \sum_{j_1, j_3 \in \{\pm1\}}  2j_1 \int_{\R}  i \rho \partial^{\ell}_x (\tilde{\psi} R_{j_3})\, N_{-j_1j_1}(\psi, \partial^{\ell}_x R_{j_1})\,dx\\[1mm]
\qquad&+\, \veps^{2}\,\mathcal{O}(\mathcal{E}_s+\veps^{3/2} \mathcal{E}^{3/2}_s)\\[2mm] 
\qquad=\;& -\veps^{2} \sum_{j_1, j_3 \in \{\pm1\}} j_1 \int_{\R}  \psi\, \tilde{\tilde{\psi}}\, \partial^{\ell}_x R_{j_1}\, \partial^{\ell+1}_x R_{j_3}\,dx\\[2mm]
\qquad&+\, \veps^{2}\,\mathcal{O}(\mathcal{E}_s+\veps^{3/2} \mathcal{E}^{3/2}_s)\,,
\end{align*}
and because of \eqref{dx1}, \eqref{dx2}, and \eqref{dxdt} we arrive at
\begin{align*}
I_3+I_4
=\;& \veps^{2}\, \int_{\R} \psi\, \tilde{\tilde{\psi}}\,   \partial^{\ell}_x (R_1+R_{-1})\,
\partial^{\ell+1}_x(R_1-R_{-1})\,dx\\[2mm]
&+\, \veps^{2}\,\mathcal{O}(\mathcal{E}_s+\veps^{3/2} \mathcal{E}^{3/2}_s)\\[2mm]
=\;& \veps^{2}\, \int_{\R}  \psi\, \tilde{\tilde{\psi}}\, \partial^{\ell}_x (R_1+R_{-1})\,
\partial_t \partial^{\ell}_x(R_1+R_{-1})\,dx\\[2mm]
&+\, \veps^{2}\,\mathcal{O}(\mathcal{E}_s+\veps^{3/2} \mathcal{E}^{3/2}_s)\\[2mm]
=\;& \frac{1}{2}\,\veps^{2}\,\partial_t \int_{\R} \psi\, (\psi + \veps^{3/2} (R_{1}+R_{-1}))\, (\partial^{\ell}_x (R_1+R_{-1}))^{2}\,dx \\[2mm]
&+\, \veps^{2}\, \mathcal{O}(\mathcal{E}_s+\veps^{3/2} \mathcal{E}^{3/2}_s)\,.
\end{align*}

Finally, we examine $I_5$. Using once more \eqref{dx1}, \eqref{dx2}, and \eqref{dxdt} yields
\begin{align*}
I_5
=\;& \frac{1}{2}\, \veps^{5/2} \sum_{j_1, j_2, j_3 \in \{\pm1\}} j_1 \int_{\R}  \partial^{\ell}_x R_{j_1}\, \partial^{\ell+1}_x (R_{j_3} R_{j_2})\,dx\\[2mm]
&+\, \veps^{5/2}\,\mathcal{O}(\mathcal{E}^{3/2}_s)\\[2mm] 
=\;& \veps^{5/2} \sum_{j_1, j_2,j_3 \in \{\pm1\}} j_1 \int_{\R} R_{j_3}\, \partial^{\ell}_x R_{j_1}\, \partial^{\ell+1}_x R_{j_2} \,dx\\[2mm]
&+\, \veps^{5/2}\,\mathcal{O}(\mathcal{E}^{3/2}_s)\\[2mm] 
=\;& -\veps^{5/2} \int_{\R}  (R_1+ R_{-1})\, \partial^{\ell}_x (R_{1}+R_{-1})\, \partial^{\ell+1}_x (R_1-R_{-1}) \,dx\\[2mm]
&+\, \veps^{5/2}\,\mathcal{O}(\mathcal{E}^{3/2}_s)\\[2mm]
=\;& - \veps^{5/2} \int_{\R}  (R_1+ R_{-1})\, \partial^{\ell}_x (R_{1}+R_{-1})\, \partial_t \partial^{\ell}_x (R_1+R_{-1}) \,dx\\[2mm]
&+\, \veps^{5/2}\,\mathcal{O}(\mathcal{E}^{3/2}_s)\\[2mm] 
=\;& -\frac{1}{2}\, \veps^{5/2} \, \partial_t \int_{\R} (R_{1}+R_{-1})\, (\partial^{\ell}_x (R_1+R_{-1}))^{2}\,dx \\[2mm]
\quad&+\, \veps^{5/2}\,\mathcal{O}(\mathcal{E}^{3/2}_s)\,.
\end{align*}
\medskip

Hence, we define the modified energy 
\begin{equation}
\tilde{\mathcal{E}_s} = \mathcal{E}_s + \frac{1}{2} \, \veps^{2}\, \sum_{\ell=1}^s h_l\,,
\end{equation}
with
\begin{equation*}
h_l = \int_{\R} \big(((2\ell+1) \partial^{2}_x \psi -\psi)(\psi + \veps^{3/2} (R_{1}+R_{-1}))+\veps^{1/2}(R_{1}+R_{-1})\big)
(\partial^{\ell}_x (R_1+R_{-1}))^{2}\,dx\,,
\end{equation*}
to obtain
\begin{equation}
\partial_t \tilde{\mathcal{E}_s}\, \lesssim\, \veps^{2} (\tilde{\mathcal{E}_s}+\veps^{1/2}\tilde{\mathcal{E}}^{3/2}_s+1)\,.
\end{equation}

Consequently, Gronwall's inequality yields for sufficiently small $\veps>0$ the $\mathcal{O}(1)$-bounded\-ness of $\tilde{\mathcal{E}_s}$ for all $t\in[0,T_0/\veps^{2}]$. Because of $\|R_{1}+R_{-1}\|_{H^{s}} \lesssim \sqrt{\tilde{\mathcal{E}_s}}$ for sufficiently small $\veps>0$ and estimate \eqref{RES2} Theorem \ref{mainresult} follows.
\qed
\medskip

\bibliographystyle{alpha}

\end{document}